\documentstyle[12pt]{article}

\begin{document}
\begin{center}
{\bf  Towards a Schubert Calculus for Maps from a Riemann Surface
to a Grassmannian }

\medskip

by

\medskip

Aaron Bertram\footnote{Partially supported by NSF grant DMS-9218215}

\end{center}

\bigskip

{\bf 0. Introduction.}
The notion of a quantum cohomology theory for complex projective varieties was
introduced
by the physicist C. Vafa (cf. \cite{V}), and has recently received a more
mathematical
treatment by Witten (cf. \cite{W}) and others (see \cite {KM} and \cite{RT}).
Given such a variety $X$ with polarization ${\cal O}_X(1)$, the idea is
to change the multiplicative structure of $\mbox{H}^*(X,{\bf C})$ by carrying
out the
intersection of cohomology classes not on $X$ itself, but rather after
pulling them back to the scheme $Mor_d(C,X)$ parametrizing the
morphisms of degree $d$ from a fixed Riemann surface $C$ to $X$.
When the right number of cohomology classes are intersected on $Mor_d(C,X)$
to produce a top codimensional cohomology class on  $Mor_d(C,X)$, then the
degree of such a class is referred to as the ``k point correlation function''
in the
literature, but following \cite{BDW},
we will simply call it the Gromov invariant of the cohomology classes.

\medskip

Although the objects in question, as described here, are all algebro-geometric,
it is not readily
apparent how to make sense of the Gromov invariants in terms of algebraic
geometry.
For example, the scheme $Mor_d(C,X)$ often does not have the expected
dimension,  is essentially never projective, and tends to have
singularities. All these features evidently pose problems when one attempts to
intersect cohomology classes.

\medskip

The purpose of this paper is first to show that if $X = G = G(r,k)$, the
Grassmannian
of complex $r$-planes in ${\bf C}^k$, then for all Riemann surfaces and
sufficiently large $d$,  the Gromov invariants for the special Schubert cycles
can be rigorously defined, and are realized as an intersection of chern classes
on a projective scheme. (One probably ought to call this ``quantum
Schubert calculus''.) It is almost immediate when one defines the
invariants in this way that they do not depend upon the choice of complex
structure on the Riemann surface of genus $g$. Another advantage of
realizing the Gromov invariants in this way comes from the general
principle that intersection numbers
tend to be more easily computed when they are defined algebraically.

This principle was demonstrated by Daskalopoulos, Wentworth
and the author in \cite{BDW}, where it was shown that in
case $r=2$ and $C$ is an elliptic curve, the Gromov invariants for
Grassmannians
agree with a remarkable conjectural formula of Vafa and Intriligator
(\cite{I}),
generalizing the classical Schubert calculus on the Grassmannian.

The second aim of this paper is to rigorously prove  an induction on the genus.
That is, an explicit relation is derived between the Gromov invariants
associated to Riemann surfaces of genus $g$ and $g-1$. The induction,
together with the results of \cite{BDW}, enables us to calculate all the
invariants when $r=2$,
and again the result agrees with the generalized Schubert calculus formula.

\medskip

The results of this paper mostly follow from a  detailed analysis of the
Grothendieck
quot scheme parametrizing quotients of a trivial vector bundle on an algebraic
curve $C$. These schemes, as was shown in \cite{BDW}, provide
compactifications
of the scheme of morphisms from $C$ to a Grassmann variety. Moreover, the pull
backs
of the special Schubert cycles
to $Mor_d(C,G)$ extend to chern classes on the quot scheme, which may therefore
be intersected in spite of the presence of singularities. It is shown that
representatives of the  pulled-back
cycles on $Mor_d(C,G)$ may be chosen in sufficiently general position so that
the actual intersection number is defined, and agrees with the
intersection of the corresponding chern classes on the quot
scheme. The induction on the genus is
obtained by considering the family of quot schemes associated to a family of
smooth
curves degenerating to an irreducible curve with one node.

\bigskip

\noindent {\bf Acknowledgements:} I would like to thank Richard Wentworth for
the many discussions we had on the topic of Gromov invariants before
and during the preparation of this paper. I would also like to
thank J\'anos Koll\'ar and Yongbin Ruan for their useful suggestions.

\newpage

\noindent {\bf 1. Intersections on the Space of Maps.} If $C$ is a Riemann
surface of
genus $g$ and $X$ is a projective variety equipped with an ample line bundle
${\cal O}_X(1)$, then the moduli space $Mor_d(C,X)$ parametrizing the morphisms
$f:C\rightarrow X$
of degree $d$ is quasiprojective and admits a universal evaluation map:
$$ev:C\times Mor_d(C,X) \rightarrow X$$
Given a subvariety $Y \subset X$ and a point $p\in C$, one gets an induced
subscheme of
$Mor_d(C,X)$ by pulling back and intersecting:
$$W_d(p,Y) := ev^{-1}(Y) \cap \left( \{p\}\times Mor_d(C,X) \right)$$
where $\{p\}\times Mor_d(C,X)$ is identified with $Mor_d(C,X)$.

\medskip

If $Y_1,...,Y_N \subset X$ are subvarieties of codimension $c_1,...,c_N$,
chosen so that
$\sum_{i=1}^N c_i = \mbox{dim}(Mor_d(C,X))$, then one expects the intersection:
$\cap _{i=1}^N W(p_i,Y_i)$ to consist of distinct reduced points.

\medskip

\noindent {\bf (Imprecise) Definition:} The Gromov invariant associated to
the subvarieties $Y_1,...,Y_N$ above is the number of points in the
intersection
of the $W_d(p_i,Y_i) \subset Mor_d(C,X)$, assuming that
$Mor_d(C,X)$ is of the
expected dimension, and that there is no extra ``intersection at infinity''.

\medskip

The Gromov invariant ought not to depend upon the choice of distinct points
$p_i \in C$.
Moreover, it may be extended by linearity to a function of cycles $Z_1,...,Z_N$
on $X$, and it ought in fact to be an invariant attached to the corresponding
cohomology classes $[Z_i] \in H^{c_i}(X,{\bf C})$.

\medskip

Our object in this section is to make rigorous sense of this definition
in case $X$ is a Grassmann variety, and to prove at least the independence of
the choice
of points.

\medskip

Suppose now that $G = G(r,k)$, the Grassmannian of $r$-planes in $V \cong {\bf
C}^k$. Let
$$0 \rightarrow S \rightarrow V\otimes {\cal O}_G \rightarrow Q \rightarrow 0$$
be the universal exact sequence of vector bundles on $G$. The polarization
${\cal O}_G(1)$  is the determinant $\wedge ^r(S^*)$, which determines the
Pl\"ucker
embedding.

\medskip

It was shown in \cite{BDW} that the dimension
of $Mor_d(C,G)$ has the expected value $kd - r(k-r)(g-1)$ if $d >> 0$, and
moreover that
in this case $Mor_d(C,G)$ is irreducible and generically reduced.

\medskip

We regard $G = G(r,k)$ as a homogeneous space for the group $GL(V)$.
If $Y\subset G(r,k)$
is any subvariety, let $gY$ be the translate of $Y$ by the element $g\in
GL(V)$.

\medskip

\noindent {\bf Lemma 1.1:} Suppose $Y\subset G(r,k)$ is an irreducible
subvariety
of codimension $c$, and suppose $Z\subset Mor_d(C,G)$ is an
irreducible subscheme. Then for any $p\in C$ and a general translate $g$,
the intersection $Z \cap W_d(p,gY)$ is either empty or has codimension $c$ in
$Z$.

\medskip

{\bf Proof:} Let $ev_p$ be the restriction of $ev$ to $\{p\}\times Mor_d(C,G)$,
and
consider the reduced image $T = ev_p(Z) \subset G(r,k)$. Under
a general translation, $gY$ intersects $T$ in codimension $c$ or the empty set.
More generally, $gY$ intersects each locus in $T$ over which $ev_p$ has
constant
fiber dimension in codimension $c$ or the empty set, so
$W_d(p,gY) \cap Z  = ev_p^{-1}(gY) \cap Z$ is of codimension $c$ in $Z$
or else empty.

\medskip

Suppose $Y_1,Y_2,...,Y_N$ are irreducible subvarieties
of $G(r,k)$ such that the codimension of $Y_i$ is $c_i$ and suppose that
$p_1,...,p_N$ are (not necessarily distinct) points of $C$.

\medskip

\noindent {\bf Corollary 1.2:}
If $\sum_{i=1}^N c_i > $dim$(Mor_d(C,G))$, then for general elements
$g_1,...,g_N \in GL(V)$, we have
$W_d(p_1,g_1Y_1)\cap ... \cap W_d(p_N,g_NY_N) = \emptyset.$

\medskip

{\bf Proof:} Immediate from the Lemma.

\medskip

Suppose in addition that $Mor_d(C,G)$ is irreducible and generically reduced.
Then:

\medskip

\noindent {\bf Corollary 1.3:} If $\sum_{i=1}^N c_i = $dim$(Mor_d(C,G))$, then
for general $g_1,...,g_N \in GL(V)$, the intersection $\cap_{i=1}^N
W_d(p_i,g_iY_i)$ is
finite, and is contained in the smooth locus.

\medskip

{\bf Proof:} Again, immediate from the Lemma.

\medskip

As a consequence of Corollary 1.3, one can count the
multiplicities at each of the finite points of the intersection
$\cap_{i=1}^N W_d(p_i,g_iY_i)$, and thus a well-defined intersection
number is obtained on $Mor_d(C,G)$.

\bigskip

One can now define the Gromov invariant associated to any subvarieties
$Y_1,...,Y_N$
as the intersection number in Corollary 1.3 for general choices of $p_1,...,p_N
\in C$
and $g_1,...,g_N\in G$. However, for special choices of $Y_1,...,Y_N$, the
intersection numbers may be realized as intersections of chern classes on a
projective scheme.

\medskip

\noindent {\bf Definition:} If $W\subset V^*$ has dimension $n \le r$, and
$Y\subset G(r,k)$ is the degeneracy locus
of the map $W\otimes {\cal O}_G \rightarrow S^*$, then following the literature
(see e.g. \cite{GH}), we will call $Y$ a special Schubert
subvariety of $G$.

\medskip

The special subvarieties of the Grassmannian are straightforward
generalizations of the
hyperplanes in projective space.  The special subvariety $Y\subset G$ is always
irreducible (though not necessarily smooth), and $Y$
represents the $r+1-n$th chern class
of $S^*$ (see \cite{GH}). Any translate $gY$ by an element of $g\in GL(V)$
is simply the degeneracy locus
associated to the translate $gW$. Finally, the special subvarieties generate
the
cohomology ring of $G$, as we will remind the reader in \S 3.

\medskip

\noindent {\bf Theorem 1.4:} If $d >>0$ and $Y_i, i =1,...,N$ are
special subvarieties of
$G(r,k)$ of codimension $c_i$ satisfying $\sum_{i=1}^N c_i = dk-(g-1)r(k-r)$,
then
the intersection number associated to $\cap _{i=1}^N W_d(p_i,g_iY_i)$ is
independent of the
choice of general $g_i\in GL(V)$, and {\it distinct} points $p_i\in C$.

\medskip

\noindent The plan of the proof of Theorem 1.4 is as follows:

\medskip

{\bf Step 1:} To find a natural (irreducible, generically reduced) projective
scheme
compactifying the space of morphisms $Mor_d(C,G)$, on which each subscheme
$W_d(p_i,g_iY_i) \subset Mor_d(C,G)$ extends
to a projective subscheme representing a chern class.

{\bf Step 2:} To show that for general $g_i$ and distinct $p_i$, the
intersection
of the $N$ extensions of the $W_d(p_i,g_iY_i)$ is contained in  $Mor_d(C,G)$,
that is,
there is no intersection at infinity.
Thus the intersection number of the theorem is computable as an intersection of
chern classes on the
compactification. Changing $g_i$  simply changes representatives of the chern
classes,
while changing the $p_i$ deforms them, and the theorem follows.

\bigskip

Following \cite{BDW}, the compactification of $Mor_d(C,G)$ we use is
Grothendieck's quot scheme. Namely, let $Quot_d(C,G)$ be the
projective scheme representing the functor:

\medskip

$F(S) = \{$ {\it flat families of quotients
of $V\otimes {\cal O}_{C\times S}$ of relative}

\hskip 1in {\it Hilbert polynomial $d-(k-r)(g-1)$ over $S$.}$\}$

\medskip

Let
$$0 \rightarrow {\cal E} \rightarrow V\otimes {\cal O}_{C\times Quot}
\rightarrow {\cal Q} \rightarrow 0$$
be the universal exact sequence on $C\times Quot_d(C,G)$. (See \cite{G} for the
details.)

While the universal quotient, ${\cal Q}$, is a rather unpleasant sheaf, the
kernel ${\cal E}$
always locally free of rank $r$ and has degree $-d$ on the fibers over
$Quot_d(C,G)$.
So  $Quot_d(C,G)$ may be thought of
dually as parametrizing injective maps $E\rightarrow V\otimes {\cal O}_C$ from
vector bundles $E$ of degree $-d$ and rank $r$.

Of course, the morphisms to the Grassmannian are special cases of this
where the map $E\rightarrow V\otimes {\cal O}_C$ is injective on each fiber,
and one easily sees that this determines $Mor_d(C,G)$ as an open subscheme
of $Quot_d(C,G)$. It was shown in \cite{BDW} moreover that for $d >>0$,
the quot scheme is irreducible and generically reduced, so $Mor_d(C,G)$ is
dense.

Now suppose that $W\subset V^*$ is a plane of dimension $n$ and
$Y \subset G(r,k)$ is the associated special subvariety. Given $p\in C$, we can
define a corresponding degeneracy locus inside $Quot_d(C,G)$ as follows.
Restrict the universal map $V^*\otimes {\cal O}_{C\times Quot} \rightarrow
{\cal E}^*$  to $p\times Quot_d(C,G)
\cong Quot_d(C,G)$. Let ${\cal E}_p$ be the restriction of ${\cal E}$ to
$p\times
Quot_d(C,G)$, and let $V_d(p,Y)$ be the degeneracy locus of the resulting
map $W\otimes {\cal O}_{Quot} \rightarrow {\cal E}_p^*$.

It is immediate from the definitions that $V_d(p,Y) \cap Mor_d(C,G) =
W_d(p,Y)$. Thus Step 1 is complete
once we show that every component of $V_d(p,Y)$ has codimension $r+1-n$ in
$Quot_d(C,G)$, since in that case $V_d(p,Y)$ will represent the $r+1-n$th chern
class of ${\cal E}_p^*$.

\bigskip

Let $Pl: G(r,k) \rightarrow {\bf P}^n$ be the Pl\"ucker embedding of the
Grassmannian
(so $n = \left( k \atop r \right) - 1$). This induces an embedding of moduli
spaces of
morphisms $Pl: Mor_d(C,G) \hookrightarrow Mor_d(C,{\bf P}^n)$, and this extends
to a morphism of quot schemes $Pl:Quot_d(C,G) \rightarrow Quot_d(C,{\bf P}^n)$
by
sending an injective map $E\rightarrow {\cal O}_C^k$ to the (injective)
deteminant map
$\wedge ^rE \rightarrow \wedge ^r {\cal O}_C^k$. One readily checks that this
induces
a map of quot functors, hence of the quot schemes.

The boundary $Quot_d(C,{\bf P}^n) - Mor_d(C,{\bf P}^n)$ decomposes as a
disjoint
union of locally closed subschemes $\sqcup _{0 < m \le d} C_m\times
Mor_{d-m}(C,{\bf P}^n)$,
where $C_m$ is the $m$th symmetric product of $C$.
This is because any map $L\rightarrow {\cal O}_C^k$ can be recovered from its
divisor $D$ of zeroes and the induced map $L(D) \rightarrow {\cal O}_C^k$.

As a result, the quot scheme decomposes as
$$Quot_d(C,G) = Mor_d(C,G) \sqcup _{0 < m \le d} B_m$$
where $B_m = Pl^{-1}(C_m\times Mor_{d-m}(C,G))$. One checks that the fibers of
the
morphisms $Pl: B_m \rightarrow C_m\times Mor_{d-m}(C,G)$ all have dimension
$(r-1)m$. In
fact, the schemes $B_m$ are all fiber bundles over $Mor_{d-m}(C,G)$ with fiber
over $E \hookrightarrow V\otimes {\cal O}_C$ isomorphic to the quot scheme
$Quot(E,b)$ of length $b$ quotients of $E$.

\medskip

Each intersection $V_d(p,Y) \cap B_m$ is contained in the union of two schemes,
namely $Pl^{-1}(C_m\times W_{d-m}(p,Y))$
and $Pl^{-1}(p+C_{m-1} \times Mor_{d-m}(C,G))$. The first scheme always has
codimension
$r+1-n$ in $B_m$ by Lemma 1.1, and the intersection of $V_d(p,Y)$ with the
second has codimension $r-n$ in the fibers of $Pl$, hence codimension $r-n+1$
in
$B_m$. This proves Step 1, in particular.

\medskip

But this setup also gives Step 2. Suppose that $W_d(p_i,g_iY_i), i = 1,...,N$
are chosen as in
the Theorem, and consider $\cap_{i=1}^N V_d(p_i,g_iY_i) \cap B_m$.
In the analysis of the previous paragraph, for at most $m$ of the points $p_i$
may the
intersection $V_d(p_i,g_iY_i)\cap B_m$ lie in the scheme
$Pl^{-1}(p_i+C_{m-1}\times
Mor_{d-m}(C,G))$ {\it because the points are distinct}. The rest must lie in
$Pl^{-1}(C_m\times W_{d-m}(p_i,g_iY_i))$. But this means, after reordering the
points,
that the intersection $\cap_{i=1}^NV_d(p_i,g_iY_i)\cap B_m$ is contained in
$Pl^{-1}(C_m\times \cap _{i=1}^{N-m}W_{d-m}(p_i,g_iY_i))$.

Suppose that $d_1$ is chosen so that $d > d_1$ implies that
$\mbox{dim}(Quot_d(C,G)) = kd -
r(k-r)(g-1)$. (Such a $d_1$ was proven to exist in \cite{BDW}). Let
$D = kd_1$. Then Step 2 follows from:

\medskip

\noindent {\bf Claim:} If $d > d_1+D$ and $g_1Y_1,...,g_NY_N$ are chosen as in
the Theorem,
then for all $0 < m \le N$, $\cap_{i=1}^{N-m} W_{d-m}(p_i,g_iY_i) = \emptyset$.

\medskip

{\bf Proof:} For $m \le D$, this follows immediately from Corollary 1.2 since
$$\begin{array}{lcl} \sum_{i=1}^{N-m} \mbox{codim}_G(Y_i) &
\ge & \mbox{dim}(Quot_d(C,G)) - mr \\
& > & \mbox{dim}(Quot_{d-m}(C,G))\end{array}$$

If $m > D$, only the last equality may fail to hold, since $Quot_{d-m}(C,G)$
may have
larger than the expected dimension, but in this case,
$$\begin{array}{lclcl} \mbox{dim}(Quot_d(C,G)) - mr &
> & D - r(k-r)(g-1) \\
& > & \mbox{dim}(Quot_{d_1}(C,G)) & \ge &
\mbox{dim}(Quot_{d-m}(C,G))\end{array}$$
and again Corollary 1.2 gives the claim.

\bigskip

Thus, the intersection number for $\cap _{i=1}^N W_d(p_i,g_iY_i)$ of Theorem
1.4
may be thought of as an invariant associated to the subvarieties $Y_1,...,Y_N
\subset G$,
or equivalently to the chern classes of $S^*$ represented by the $Y_i$.

\bigskip

\noindent{\bf Proposition 1.5:} The intersection number of Theorem 1.4
does not depend upon the choice of smooth curve of genus $g$.

\medskip

{\bf Proof:} In light of the irreducibility of the moduli space of smooth
curves of genus
$g$, we have only to prove the following. For each family ${\cal C}$ of smooth
curves over a
connected base $B$, the intersection numbers in Theorem 1.4 are the same for
all fibers
${\cal C}(b)$.

Let ${\cal Q}_d$ be the relative quot scheme over the base $B$, with fibers
isomorphic
to $Quot_d({\cal C}(b),G)$. The fact that the quot schemes are all of the
expected dimension
($d$ is assumed to be large) implies that the map ${\cal Q}_d \rightarrow B$ is
a complete intersection morphism, hence the family of quot schemes is flat over
$B$
(see e.g.\cite{K}).

Let ${\cal U} \hookrightarrow V\otimes {\cal O}_{{\cal C}\times_B {\cal Q}_d}$
be the
universal subbundle. After restriction and base change, we can find a
section $\sigma$ of ${\cal C}$ near each point $b\in B$, and
the restriction of ${\cal U}^*$ to $\sigma \times_B {\cal Q}_d$
gives a varying (flat) family of vector bundles over a neighborhood of $b$. The
Proposition now follows
from:

\medskip

\noindent {\bf Lemma 1.6:} If $f: {\cal X} \rightarrow {\cal Y}$ is a
projective, flat
morphism of relative dimension $n$ over an irreducible base, and if
${\cal E}$ is a vector bundle of rank $r$
on ${\cal X}$, then for any polynomial $P(X_1,...,X_r)$ of weighted degree $n$,
the
intersection number $\int_{{\cal X}_y}P(c_1({\cal E}_y),...,c_r({\cal E}_y))$
is
independent of the point $y\in {\cal Y}$.

\medskip

{\bf Proof:} Let ${\cal O}_{\cal X}(1)$ be a relatively ample line bundle on
${\cal X}$. If $L_1,...,L_k$ are line bundles on ${\cal X}$ and $\sum_{i=1}^k
a_i = n$,
then it follows from the invariance of the Hilbert polynomials:
$\chi ({\cal X}_y,\otimes {L_i}^{ma_i} (n))$ that the intersection numbers
$c_1(L_1)^{a_1}...c_1(L_k)^{a_k}.[{\cal X}_y]$  are independent of $y\in {\cal
Y}$.
If we apply this to the appropriate intersections of pullbacks of ${\cal
O}_{\bf P}(1)$'s
to fiber products of projective bundles ${\bf P}({\cal E})\times ... \times
{\bf P}({\cal E})$ over ${\cal Y}$,
then we get invariance of the intersections of segre classes, hence of the
chern classes
as well.

\bigskip

Suppose that $X_1,...,X_r$ are weighted variables such that the weight of $X_i$
is $i$.
If $P(X_1,...,X_r)$ is a homogeneous polynomial of weighted degree $kd -
r(k-r)(g-1)$,
then for $d >>0$:

\medskip

\noindent {\bf Definition:} The genus-$g$ Gromov invariant
$N_d(P(X_1,...,X_r),g)$ is
the intersection number of Theorem 1.4, extended by linearity from the
monomials,
where $X_i$ are replaced by special subvarieties of codimension $i$.

\bigskip

\noindent {\bf Warning:} The Gromov invariants do not respect intersections!
Namely, as noted following Corollary 1.3, it is possible to attach
a Gromov invariant to any subvarieties $Y_1,...,Y_N$ of $G$. In general,
if $Y_1$ and $Y_2$, for example, intersect transversally, then
the Gromov invariant for $Y_1\cap Y_2,Y_3,...,Y_N$ will not coincide with the
Gromov
invariant for $Y_1,...,Y_N$.

\bigskip

\noindent {\bf Remark:} It was shown in \cite{BDW} that one can define a
consistent
Gromov invariant for all degrees by downward induction from the observation
that
$$N_d(P(X_1,...,X_r),g) = N_{d-r}(X_r^kP(X_1,...,X_r),g)$$
for sufficiently large $d$.

\bigskip

\noindent {\bf Proposition 1.7:} If $C_0$ is an irreducible curve with a single
node $\nu$
and arithmetic genus $g$, then the
Gromov invariants for $C_0$ and special subvarieties of $G$
agree with the genus-$g$ Gromov invariants.

\medskip

{\bf Proof:} This time, we use a family of curves
${\cal C} \rightarrow B$
smoothing the nodal curve $C_0$, and use the invariance property of Lemma 1.6
to
conclude that if ${\cal E} \hookrightarrow V\otimes {\cal O}_{C_0\times Quot}$
is the universal
(torsion-free) subsheaf and $p\in C_0 - \nu$,
then the intersection numbers of chern classes of the (locally free) sheaf
${\cal E}_p$ are the same
as the corresponding intersection numbers in the smooth case.

\medskip

The argument proceeds just as in the proofs of Theorem 1.4 and Proposition 1.5.
That is, we need to show that:

(1) For $d >>0$, the quot scheme $Quot_d(C_0,G)$ is irreducible and generically
reduced
of the expected dimension $dk-r(k-r)(g-1)$.

\smallskip

(2) For $d >>0$, the degeneracy loci $V_d(p_i,g_iY_i)$ represent chern classes
of ${\cal E}_p$, and the intersections $\cap_{i=1}^N V_d(p_i,g_iY_i)$ and
$\cap_{i=1}^N W_d(p_i,g_iY_i)$ are the same.

\medskip

We need to use the following result on the structure of torsion-free sheaves on
$C_0$,
which may be found, for example, in \cite{S}. Namely, if $E$ is a torsion-free
sheaf,
then the cokernel, $T_\nu$, of the double-dual exact sequence:
$$0 \rightarrow E \rightarrow E^{**} \rightarrow T_\nu \rightarrow 0$$
is a sheaf (necessarily of rank $< r$) supported on the node $\nu$.

In other words,
the torsion-free sheaves of rank $r$ and degree $-d$ are identified with the
locally free sheaves $F$ of rank $r$ and degree $-d+m$, together with a
quotient
of the fiber $F(\nu)$ of rank $m \le r$.

\medskip

This implies that the quot scheme $Quot_d(C_0,G)$ decomposes as a disjoint
union of
locally closed subschemes:
$$Quot_d(C_0,G) = Mor_d(C_0,G) \sqcup _{m=1}^d B_m \sqcup_{m=1}^r A_m$$
where $Pl: B_m \rightarrow Mor_{d-m}(C_0,G)$ is smooth, of relative dimension
$mr$
as before, and
the new pieces $A_m$ are $G(m,r)$-bundles over $Quot^{lf}_{d-m}(C_0,G)
:= Mor_{d-m}(C_0,G) \sqcup _{i=1}^{d-m} B_i$, the subscheme
of $Quot_{d-m}(C_0,G)$ parametrizing locally-free subsheaves
$E\hookrightarrow V\otimes {\cal O}_{C_0}$.

\medskip

The same argument as in the proof of Theorem 1.4 will give (1) and (2) once we
show:

\medskip

\noindent {\bf Claim:} $Mor_d(C_0,G)$ is of the expected dimension for $d >>
0$.

\medskip

{\bf Proof:} Let $f: \tilde C_0 \rightarrow C_0$ be the normalization of $C_0$.
It suffices to show
that the image of the obvious embedding:
$\iota: Mor_d(C_0,G) \hookrightarrow Mor_d(\tilde C_0,G)$ has
codimension  $r(k-r)$ for $d >> 0$.

The codimension of the locus $\{ E \rightarrow V\otimes {\cal O}_{\tilde C_0} |
H^1(\tilde C_0,
E^*(-p-q)) \ne 0 \}$ in $Mor_d(\tilde C_0,G)$ grows with $d$. This follows, for
example, from
the proof of Theorem 4.28 in \cite {BDW}. In particular, for large $d$, it
exceeds
$r(k-r)$, and thus this locus may be ignored.

On the other hand, if $H^1(\tilde C_0,E^*(-p-q)) = 0$, then the natural map
$\mbox{Hom}(E,V\otimes {\cal O}_{\tilde C_0}) \rightarrow \mbox{Hom}(E_p,V)
\oplus \mbox{Hom}(E_q,V)$
is surjective, and the claim follows by pulling back the diagonal from the
product of Grassmannians.

\bigskip

{\bf 2. Induction on the Genus.} The universal exact sequence implies that
the chern classes of the universal quotient bundle
$Q$ on $G$ are all expressible in terms of $c_1(S^*),...,c_r(S^*)$, and the
chern classes of the tangent bundle $TG = S^* \otimes Q$ are therefore
expressed
in terms of $c_1(S^*),...,c_r(S^*)$ by the standard formulas for the chern
classes
of a tensor product. (See \S 3 for more details.) In particular, we get in this
way
an ``euler'' polynomial $e(X_1,...,X_r)$ of weighted degree $r(k-r)$,
with the property that when evaluated at the chern
classes of $S^*$, $e$ produces the euler class, that is, the top chern class
of $TG$. The point of this section is to prove the following induction formula:

\bigskip

\noindent {\bf Theorem 2.1:} The Gromov invariants for genus $g$ and $g-1$ are
related by:
$$N_d(P(X_1,...,X_r),g) = N_d(e(X_1,...,X_r)P(X_1,...,X_r),g-1)$$

\bigskip

The idea of the proof is as follows. We have already seen in Proposition 1.7
that
the intersection $\cap_{i=1}^NW_d(p_i,g_iY_i)$ for an irreducible curve $C_0$
with one node and
arithmetic genus $g$ computes the Gromov invariant for genus $g$. On the other
hand,
the scheme $Mor_d(C_0,G)$ embeds in $Mor_d(\tilde C_0,G)$ as the subscheme
parametrizing morphisms which send $p$ and $q$ to the same point. If the other
points
$p_i\in C_0$ are identified with their preimages in $\tilde C_0$, then
the schemes $W_d(p_i,g_iY_i) \subset Mor_d(C_0,G)$ are just the intersections
of the corresponding subschemes in $Mor_d(\tilde C_0,G)$ with the image of
$Mor_d(C_0,G)$.

The proof of Theorem 2.1 consists in showing that $Mor_d(C_0,G)$ extends to
a subscheme of $Quot_d(\tilde C_0,G)$ which represents the pullback from $G$
of the euler polynomial in $c_i(S^*)$, and that there are no
``intersections at infinity'' when this scheme is intersected with the
$V_d(p_i,g_iY_i)$.

\bigskip

The first problem we encounter is the fact that $ev_p^*(Q)$ does not extend as
a vector bundle to the quot scheme.

\bigskip

\noindent {\bf Definitions:} (a) Let $G^* = G(k-r,k)$ be the Grassmannian of
$k-r$-dimensional
subspaces of $V^*$ and let ${\cal F} \hookrightarrow V^*\otimes {\cal O}$ be
the universal
subbundle on $C\times Quot_d(C,G^*)$.

(b) Let $\Gamma \subset Quot_d(C,G) \times Quot_d(C,G^*)$ be the closure of the
image of
$Mor_d(C,G) = Mor_d(C,G^*)$ via the two embeddings. Let the two
projections from $\Gamma$ be named $\pi: \Gamma \rightarrow Quot_d(C,G)$
and $\pi^*:\Gamma \rightarrow Quot_d(C,G^*)$.

(c) Let $Z_{p,q}$ be the zero scheme in $\Gamma$ of the canonical
map ${\cal E}_p \rightarrow {\cal F}_q^*$
(which factors through $V\otimes {\cal O}_{C\times \Gamma}$.)

\bigskip

Note that $Mor_d(C_0,G)$ is identified with an open subscheme of $Z_{p,q}$
via $\pi^{-1}$.

\medskip

Let $M(X_1,...,X_r)$ be a monomial of weighted degree $kd-r(k-r)(g-1)$ with
corresponding
special subvarieties
$Y_1,...,Y_N$. Let $c_i = c_i({\cal E}^*_r)$ for some point $r\in \tilde C_0$.
The theorem follows immediately from:

\medskip

\noindent {\bf Lemma 2.2:} (a) The following degrees are the same:
$$M(c_1,...,c_r).[Z_{p,q}]
= e(c_1,...,c_r)M(c_1,...,c_r).[Quot_d(\tilde C_0,G)]$$

\medskip

(b) For general $p_1,...,p_N \in \tilde C_0$ and $g_1,...,g_N\in GL(V)$,
$$Z_{p,q} \cap_{i=1}^N V_d(p_i,g_iY_i) = Mor_d(C_0,G)
\cap_{i=1}^NW_d(p_i,g_iY_i)$$

\medskip

{\bf Proof of (a):}
Since $Z_{p,q}$ is the zero locus of a section $\sigma \in Hom({\cal E}_p,{\cal
F}_q^*)$,
it follows that if $Z_{p,q}$ is irreducible, of
codimension $r(k-r)$, then $Z_{p,q}$ represents the top chern
class of ${\cal E}^*_p\otimes {\cal F}^*_q$. Moreover, since ${\cal E}^*_p
\otimes {\cal F}^*_q$
coincides with the pullback $ev_p^*S^* \otimes ev_q^*Q$ from $G$ over the open
subscheme $Mor_d(\tilde C_0,G)$, any difference between $Z_{p,q}$ and the
pull-back
of the euler
polynomial would be concentrated in higher codimension than $r(k-r)$, hence
can be ignored.

\medskip

The same argument as in the proof of Proposition 1.7 shows that $Z_{p,q}$
intersects $Mor_d(\tilde C_0,G)$ in codimension $r(k-r)$.
As with Step 1 in the proof of Theorem 1.4, we next need to analyze the
intersection of
$Z_{p,q}$  with the boundary subschemes
$\pi^{-1}B_m = \pi^{-1}Pl^{-1}(C_m \times Mor_{d-m})$ for $1 \le m \le d$.
A new argument is required only for
the intersections of $Z_{p,q}$ with the open subschemes of
$Pl^{-1}(ap+bq\times Mor_{d-a-b})$ parametrizing maps $E \hookrightarrow
V\otimes {\cal O}_C$
which have rank $r-a$ at $p$ and $r-b$ at $q$.

If $F$ is a vector bundle on $\tilde C_0$ of degree $-d+a+b$ with
$\mbox{H}^1(F^*(-p-q)) = 0$, then the codimension in
$$\{ s: F \hookrightarrow V\otimes {\cal O}| rk(s_p) = rk(s_q) = r \}$$
of the set of such $s$ so that $Z_{p,q}\cap \pi^{-1}Pl^{-1}(s) \ne \emptyset$
is the same as the codimension in $G(r,k)\times G(r,k)$ of
$\{ (\Lambda,\Lambda')| \mbox{dim}(\Lambda \cap \Lambda') \ge r-a-b \}$, which
is $(k-r-a-b)(r-a-b)$. From this it follows that the codimension of
$Z_{p,q} \cap \pi^{-1}Pl^{-1}(ap+bq\times Mor_{d-a-b})$ in $\Gamma$ is
$r(k-r)+a^2+b^2$.

This argument readily adapts to show that all other intersections with boundary
components also have codimension larger than $r(k-r)$, and part (a) follows.

\medskip

{\bf Proof of (b):} As in the claim in the proof of Theorem 1.4, we may assume
that all $Mor_{d-m}(\tilde C_0,G)$ are of dimension $dk - r(k-r)(g-1) - mk$.
But then by the
analysis of the previous paragraph, the locus in $ap + bq \times
Mor_{d-a-b}(\tilde C_0,G)$ over
which $Z_{p,q}$ intersects nontrivially has codimension $r(k-r) - (a+b)k +
(a+b)^2$, and
as in the proof of Step 2, this fact, together with Corollary 1.2, implies part
(b).

\bigskip

{\bf 3. Quantum Schubert Calculus.} Both the Chow ring and the
cohomology ring of $G(r,k)$
are isomorphic to ${\bf C}[X_1,...,X_r]/I$,
where $I$ is the ideal
generated by the coefficients of $t^{k-r+1},...,t^k$ in the
formal-power-series inverse of $$P_t := 1+X_1t+X_2t^2+...+X_rt^r$$
The isomorphisms may be realized by identifying $X_i$ with the $i$th chern
class
of $S^*$. In particular, the $X_i$ are represented by the special Schubert
subvarieties.

\medskip

If $X_i$ is assigned weight $i$ and if $P(X_1,...,X_r)$ is a
polynomial of weighted degree
$r(k-r) = {\rm dim}(G(r,k))$, then the computation of
$\int_G P(c_1(S^*),...,c_r(S^*))$ is an intersection of special subvarieties,
and
is a special case of the Schubert calculus on the
Grassmannian. The classical Pieri formulas
compute these integers (see \cite{F} or \cite{GH}), but a recent formula
due to Vafa and Intriligator (\cite{I}) generalizes readily to the
quantum case.

\medskip

Let $Q_t = 1 + \sum _{m=1}^\infty y_mt^m$ be the formal-power-series
inverse of $P_t$. Then the polynomials $y_m(X_1,...,X_r)$ can
all be computed from the expansion
$$\mbox {log}(P_t) = \sum_{n=0}^\infty W_n(X_1,...,X_r)t^n.$$
Namely, differentiating both
sides with respect to $X_i$ and equating powers of $t^n$, one finds that
$y_{n-i} = \frac {\partial W_n}{\partial X_i}$ for $1 \le i \le r$.

In particular, (plus or minus) the coefficients $y_{k-r+1},...,y_k$, which
generate
all of the relations in the cohomology ring,
are the partial derivatives of a single
polynomial $W := (-1)^kW_{k+1}(X_1,...,X_r)$. It is immediate that $W$ has
weighted
degree $k+1$. Indeed, if $q_1,...,q_r$ are the chern roots of $P_t$, that is
$q_1,...,q_r$
are formal variables satisfying $P_t = \prod _{i=1}^r (1 + q_it)$, then in
terms
of these chern roots, we may write $$W = \sum_{i=1}^r \frac {q_i^{k+1}}{k+1}.$$

If we define the perturbed polynomial
$$\widetilde W := W + X_1,$$
then the zeroes
of the $X_i$-gradient $\nabla _X \widetilde W$ are distinct and reduced.
Indeed, the
zeroes of the $q_i$-gradient $\nabla_q \widetilde W$ are reduced, and occur
when
each $q_i$ is a $k$th root of $-1$. The change of variables coming from the
elementary
symmetric polynomials $X_i = \sigma _i(q_1,...,q_r)$ implies that the zeroes of
$\nabla_X\widetilde W$ are reduced and that there are $\left( k \atop r
\right)$ of them (one
for each unordered $r$-tuple of distinct $k$th roots of $-1$).

We let $$h(X_1,...,X_r) =
\mbox{det}\left( \frac {\partial ^2 \widetilde W}{\partial X_i \partial X_j}
\right)$$ be the
Hessian polynomial. Then Vafa and Intriligator tell us:

\bigskip

\noindent {\bf Proposition 3.1:} The Schubert calculus for special subvarieties
of $G(r,k)$
is computed by the following formula:
$$(*) \ \ \int_G P(c_1(S^*),...,c_r(S^*)) = (-1)^{\left( r\atop 2 \right)}
\sum_{\nabla_X\widetilde W = 0} Ph^{-1}$$
where $P(X_1,...,X_r)$ is a homogeneous polynomial of weighted degree $r(k-r)$.

\bigskip

{\bf Proof:} Since both the left and right-hand sides of $(*)$ are linear in
$P$,
and $\mbox H^{r(k-r)}(G,{\bf C})$ is one-dimensional, it suffices to show:

\medskip

(1) $(*)$ holds and is nonzero for one choice of $P$.

\medskip

(2) If $P \in I$, where $I$ is the ideal of relations among the $c_i(S^*)$,
then the right-hand
side evaluates to zero.

\bigskip

The evaluation of $h(X_1,...,X_r)$ (or $h^{-1}(X_1,...,X_r)$) on the zeroes of
$\nabla _X \widetilde W$ may be written in terms of the $q_i$'s as follows.
The determinant of the Jacobian matrix $J$ associated
to the elementary symmetric functions $\sigma_i(q_1,...,q_r)$ is the
Vandermonde determinant
$\prod_{i<j}(q_i - q_j)$, so $\nabla_X \widetilde W = (\nabla_q \widetilde
W)(J^{-1})$,
and when evaluated at the critical points $\nabla_X\widetilde W = 0$, we have
$$h(X_1,...,X_r) = \mbox{det}\left( \frac{\partial ^2 \widetilde W}{\partial
q_i \partial q_j}
\right) \mbox{det}(J^{-2}) = \frac {k^r \prod_{i=1}^r q_i^{k-1}}{(\prod_{i <
j}(q_i-q_j))^2}$$

\medskip

We verify (1) for the choice of $P(X_1,...,X_r) = X_r ^{k-r}$. The
Schubert cycle $Y \subset G(r,k)$ consisting of the zero locus of a section
${\cal O}_G \rightarrow S^*$ represents $c_r(S^*)$, and $k-r$ generally
chosen sections yield Schubert cycles intersecting in a single point,
corresponding to
the plane $W\subset V^*$ spanned by the sections.
Thus, $\int_G c_r^{k-r}(S^*)  = 1$, that is, this choice of $P$ gives
the volume form on the Grassmannian.

We compute the right-hand-side using the $q$ variables. Namely,
$$\sum_{\nabla \widetilde W = 0} X_r^{k-r} h^{-1}
  =   \sum_{\nabla \widetilde W = 0} \prod_{i=1}^r q_i^{k-r}
\frac {(\prod_{i< j}(q_i-q_j))^2}{k^r \prod_{i=1}^r q_i^{k-1}} $$
$$ = \frac {1}{r!k^r} \sum_{(-q_i)^k = -1} \prod _{i=1}^r q_i^{1-r}
(\prod_{i < j} (q_i - q_j))^2$$
Only the term $q_1^{r-1}\cdots q_r^{r-1}$ in $(\prod_{i< j}(q_i - q_j))^2$
contributes
nontrivially to the sum. This term appears with coefficient $(-1)^{\frac
{r(r-1)}2}(r!)$,
yielding $(-1)^{\frac {r(r-1)}2}$ as the sum, as desired.

\medskip

To verify (2), it is enough by linearity to show that if $P = \left( \frac
{\partial W}
{\partial X_i} \right) N$, then $\sum_{\nabla \widetilde W = 0} Ph^{-1} = 0$.
But
$\frac {\partial W}{\partial X_i} = \frac {\partial \widetilde W}{\partial
X_i}$ for $ i \ge 2$,
so the result is only nontrivial in case $i = 1$. Since
$\frac {\partial \widetilde W}{\partial X_1} = \frac {\partial W}{\partial X_1}
+ 1$,
it follows that
if $P = \left( \frac {\partial W}{\partial X_1} \right) \cdot N$, then
$\sum_{\nabla \widetilde W = 0} Ph^{-1} = - \sum_{\nabla \widetilde W =
0}Nh^{-1}$, and
when written in terms of the $q_i$'s,
this sum is easily seen to be zero.

\bigskip

\noindent {\bf Corollary 3.2:} The euler class $c_{r(k-r)}(TG)$ and the
(modified)
hessian polynomial $(-1)^{\left( r \atop 2 \right)}h(c_1(S^*),...,c_r(S^*))$
determine the
same element of $\mbox{H}^*(G)$.

\bigskip

{\bf Proof:} We have determined that the number of zeroes of $\nabla_X
\widetilde W$ is
$\left( k \atop r \right)$. On the other hand,
(23.2.1) of \cite{BT} states that the Poincar\'e series of $G(r,k)$ is:
$$(1-t^2)...(1-t^{2k})/(1-t^2)...(1-t^{2r})(1-t^2)...(1-t^{2(k-r)})$$
Thus it follows from Proposition 2.1 that:
$$(-1)^{\left( r \atop 2 \right)} \int_G h(c_1(S^*),...,c_r(S^*))
= \left( k \atop r \right) = \int_G c_{r(k-r)}(TG) $$

\bigskip

This says that the euler class, as a polynomial in the chern classes of $S^*$,
and the modified hessian are the same modulo the ideal generated by
$\frac {\partial W}{\partial X_1}, ... ,\frac {\partial W}{\partial X_r}$.
But the following improvement
will show that they are in fact the same in the quantum cohomology ring.

\bigskip

\noindent {\bf Proposition 3.3:} The euler polynomial and $(-1)^{\left( r \atop
2 \right) }h$
are the same at the zeroes of $\nabla \widetilde W$.

\medskip

{\bf Proof:} The splitting principle together with the tensor product formula
imply
that $$(**) \ e(X_1,...,X_r) = c_{r(k-r)}(S^*\otimes Q) = \prod_{i=1}^r
(q_i^{k-r} +
q_i^{k-r-1}c_1(Q) + ... + c_{k-r}(Q))$$

When evaluated at the points where $\nabla \widetilde W = 0$, we have:
$$\begin{array}{rcl} \frac 1{\prod_{j=1}^r(1-q_jt)} & = & \frac 1{c_t(S)}\\
 & = & 1 + c_1(Q)t + ... + c_{k-r}(Q)t^{k-r} + 0 + ... + 0
- t^k - c_1(Q)t^{k+1} - ...\\
& = & (1 + c_1(Q)t + ... + c_{k-r}(Q)t^{k-r})(1 + t^k)^{-1}
\end{array}$$

If we take residues of both sides at $t=q_i^{-1}$, where $q_i$ is a kth root of
$-1$,
we get:
$$\frac 1{\prod_{j\ne i} (1-q_jq_i^{-1})} =
\frac 1k(1 + c_1(Q)q_i^{-1} + ... + c_{k-r}(Q)q_i^{-(k-r)})$$

Multiplying by $kq_i^{k-r}$ and taking the product over all $i$, the right hand
side
gives the euler polynomial by $(**)$, but the left hand side gives
$(-1)^{\left( r \atop 2 \right)}$ times the formula for
the hessian in terms of the $q_i$'s derived earlier.

\bigskip

Putting together the results of this paper with \cite{BDW},
we get the following formula, first conjectured
by Vafa and Intriligator in \cite{I}:

\medskip

\noindent {\bf Theorem 3.4:} The Gromov invariants for $G = G(2,k)$
are computed by the formula:
$$N_d(P(X_1,X_2),g) = (-1)^{g-1}\sum_{\nabla \widetilde W = 0}
P(X_1,X_2)h(X_1,X_2)^{g-1}$$

\medskip

{\bf Proof:} As already mentioned, the main result of \cite{BDW} shows that
the Theorem holds for curves of genus one, and the arguments therein also
imply the genus zero case. For higher genera, we use Theorem 2.1, Proposition
3.3
and the result for genus one to get:
$$N_d(P,g) = N_d(Pe^{g-1},1) = \sum_{\nabla \widetilde W = 0}Pe^{g-1} =
(-1)^{g-1}\sum_{\nabla \widetilde W = 0}Ph^{g-1}$$
as desired.

\bigskip

The full conjecture of Vafa and Intriligator states that the analogous formula
computes the genus-$g$ Gromov invariants on any Grassmannian. Namely,

\medskip

\noindent {\bf Conjecture:} If  $P(X_1,...,X_r)$ is of
 weighted degree $kd-r(k-r)(g-1)$, then
$$ N_d(P(X_1,...,X_r),g) =
(-1)^{(g-1)\left( r \atop 2 \right)}\sum_{\nabla\widetilde W = 0} Ph^{g-1}$$

\bigskip

In light of Theorem 2.1 and Proposition 3.3, it would suffice to verify this
conjecture in the genus zero case to get the result for all genera. Since the
genus zero quot scheme is smooth and indeed has been described quite explicitly
by Str\o mme
in \cite{St}, it may be possible to verify this directly. Indeed, it was
recently
shown in \cite{R} by direct methods that the degree of the genus zero quot
scheme, namely,
$N_d(X_1^{kd+r(k-r)},0)$, agrees with the
conjectured value.

\bigskip

University of Utah, Salt Lake City, UT 84112

{\it email address:} bertram@math.utah.edu

\end{document}